\documentclass[12pt]{amsart}

\usepackage{latexsym}
\usepackage{amsmath}
\usepackage{amssymb}
\usepackage{amscd}
\usepackage{multicol}
\usepackage{multicol}
\usepackage[cmtip,matrix,arrow]{xy}
\usepackage{epsfig}

\newtheorem{theorem}{\bf Theorem}[section]

\newtheorem{lemma}[theorem]{\bf Lemma}

\newtheorem{definition-theorem}[theorem]{\bf Theorem-Definition}

\def\min{{\rm min}}

\def\P{\mathbb{P}}
\def\bR{\mathbb{R}}
\def\bC{\mathbb{C}}

\def\t{\mathfrak{t}}

\def\g{\mathfrak{g}}
 \def\k{\mathfrak{k}}
 
\def\quott{/\! /}

\def\g{\frak{g}}
\def\t{\frak{t}}

\def\u{\frak{u}}
\def\k{\frak{k}}

\def\e{e^*}

\def\O{{\mathcal O}_{\lambda}}

\title[Symplectic quotients of Schubert varieties]{On some symplectic quotients of Schubert varieties}

\author[A.-L. Mare]{Augustin-Liviu  Mare}
\address{
Department of Mathematics and Statistics\\ University of
Regina \\  Regina SK, Canada S4S
0A2}
\email{mareal@math.uregina.ca}\date{\today}

\begin{document}
\maketitle

\begin{abstract} 
Let $G/P$ be a generalized flag variety, where $G$ is a complex semisimple
connected Lie group and $P\subset G$ a parabolic subgroup.
Let also $X\subset G/P$ be a Schubert variety. 
We consider the canonical embedding of $X$ into a projective space,
which is  obtained by identifying $G/P$ with a 
coadjoint orbit of the compact Lie group $K$, where $G=K^{\bC}$.
The maximal torus $T$ of $K$ acts linearly on the projective space and
it leaves $X$ invariant: let $\Psi: X \to {\rm Lie}(T)^*$ be the restriction of the
moment map relative to the Fubini-Study symplectic form. By a theorem of Atiyah, $\Psi(X)$ is a convex polytope in
${\rm Lie}(T)^*$. In this paper we show that all pre-images $\Psi^{-1}(\mu)$,
$\mu\in \Psi(X)$, are connected subspaces of $X$. 
We then consider a  one-dimensional subtorus $S\subset T$, and the map
$f: X\to \bR$, which is the restriction of the $S$ moment map   to $X$.
We study  quotients of the form 
$f^{-1}(r)/S$,
where $r\in \bR$.
We show that under certain assumptions concerning
$X$, $S$, and $r$,  
these symplectic quotients  are (new) examples of spaces
for which the Kirwan surjectivity theorem and   Tolman and Weitsman's
presentation of the kernel of the Kirwan map hold true
(combined with a theorem of Goresky, Kottwitz, and MacPherson, these results
lead to an explicit description of the cohomology ring of the quotient). 
 The singular Schubert variety in the Grassmannian $G_2(\bC^4)$ of
 2 planes in $\bC^4$ is
 discussed in detail.

\vspace{0.5cm}

\noindent 2000 {\it Mathematics Subject Classification.} 53D20, 14L24

\end{abstract} 

\section{Introduction}
Let $K$ be a compact connected semisimple Lie group and $T\subset K$ a
maximal torus.
We also consider   the complexification $G$ of $K$ and $P\subset G$ 
a parabolic subgroup. In this paper we study Schubert varieties in the flag manifold
$G/P$. More specifically, any
such variety $X$ is $T$ invariant and admits  canonical $T$ equivariant embeddings into projective spaces
with linear $T$ actions:
we will be interested  
in  the symplectic quotients of $X$  induced by   the action of $T$.
We need to give more details  in order to be able to state the results. 
Let us denote by $R$ the 
roots of $K$ relative to $T$ and by $R^+$ the set of all positive roots with respect to a certain choice of a 
simple root system. Let also $\k$, $\g$, $\t$ denote the Lie algebras of $K$, $G$, respectively $T$. We may assume that the Lie algebra of $P$
is $\t^{\bC}\oplus \bigoplus \g_\alpha$, where $\alpha\in R^+\cup (-R^+_P)$.
Here $\g_\alpha$ is the root space of $\alpha$ and $R^+_P$ is a certain
subset of $R^+$ canonically associated to $P$ (cf. \cite[Section 23.3]{Fu-Ha}).
Let $W_P$ denote the subgroup of $W$ generated by
all reflections $s_\alpha$, with $\alpha \in R_P^+$.  The Schubert  cells in $G/P$ 
give a cell decomposition of this space. They are labeled by the quotient $W/W_P$. Namely, let us consider
the Borel subgroup $B$ of $G$ whose Lie algebra is  
$\t^\bC\oplus\bigoplus_{\alpha\in R^+}
\g_\alpha$.  The Schubert cells
are the $B$ orbits $BwP/P$, where $w\in W/W_P$. To any such $w$ corresponds the Schubert variety  
$$X(w)=\overline{BwP/P}.$$

In the  equation above, the closure  is taken relative to the differential topology of $G/P$: we will see 
that both $G/P$ and $X(w)$ are (Zariski closed) subvarieties of the same projective space. Such projective embeddings are constructed  as follows.
Let us first pick a dominant weight $\lambda\in \t^*$ such that
$$R^+_P=\{\alpha\in R^+ \ : \ \lambda(\alpha^\vee)= 0\}.$$
Here $\{\alpha^\vee \ : \alpha\in R\}$ denotes the root system dual to $R$. 
Then we have $$P=P_\lambda,$$
that is, the Lie subgroup of $G$ whose Lie algebra is
$\t^{\bC}\oplus\bigoplus_{\alpha\in R, \lambda(\alpha^\vee)\ge 0}\g_\alpha.$ 
Let $\chi_{\lambda} : P_\lambda \to \bC^{*}$ be the group homomorphism whose differential  
$d(\chi_{\lambda})_e:{\rm Lie}(P_\lambda) \to \bC$ is  the composition of $\lambda \otimes \bC$
(regarded as a $\bC$-linear function on $\t^\bC$)
with the natural projection ${\rm Lie}(P_\lambda) \to {\rm Lie}(T^{\bC})$. Consider the line bundle $L_{\lambda}$ over $G/P_\lambda$ whose total space is 
$G \times \bC/ P_\lambda$, where $h.(g,z) := (gh^{-1},\chi_{\lambda}^{-1}(h)z)$, for 
$h\in P_\lambda$, $(g,z)\in G\times \bC$. One can show that $L_{\lambda}$ is very ample. More concretely, one can embed $G/P_\lambda$ into $\P(V_{\lambda})$
as a $G$ orbit, in such a way that $L_{\lambda}$ is the restriction to $G/P_\lambda$ of the hyperplane bundle over $\P(V_{\lambda})$. Here  
$V_{\lambda}=\Gamma_{\it hol}(G/P_\lambda, L_{\lambda})^*$ is the  irreducible representation of $K$ of highest weight $\lambda$: this is the content of the Borel-Weil theorem  (see e.g. 
\cite[Section 4.12]{Du-Ko}, \cite[Chapter V]{Huc}). The action of $B$ on 
$\P(V_\lambda)$ is linear, thus the
Schubert variety $X(w)$ defined above  is a (Zariski closed) subvariety of
$\P(V_\lambda)$.  Note  that $X(w)$ is in general not smooth. The  singularities of Schubert varieties have been intensively investigated  (see for instance the survey \cite{Bi-La}). For instance, by a theorem of Ramanathan \cite{Ra},  
$X(w)$ has  rational singularities.

We equip $\P(V_\lambda)$ with the Fubini-Study symplectic form.
A moment map of the $T$ action is 
 $\Psi:\P(V_\lambda)\to \t^*$, whose component along $\xi\in\t$ is given by
 $$\Psi_\xi([v])= \frac{i}{2\pi}\frac{(\xi v)\cdot v}{v\cdot v},$$
 for all $v\in V_\lambda\setminus\{0\}$.
 Here ``$\cdot$" is a $K$-invariant Hermitean inner product on $V_\lambda$,
 which is uniquely determined up to a non-zero factor; the vector $\xi v\in V_\lambda$ arises from the infinitesimal automorphism of $V_\lambda$ induced by $\xi$.
As already mentioned, we will study symplectic quotients of $X(w)$ relative to
 $\Psi|_{X(w)}$. Of course this map is the same as the restriction of
 $\Psi|_{G/P_\lambda}$ to $X(w)$. This observation is useful to us, 
 since $\Psi|_{G/P_\lambda}$ is well understood. Namely, we identify
\begin{equation}\label{ide}G/P_\lambda=\O,\end{equation} 
 where  $\O:={\rm Ad}^*(K)\lambda$ is the coadjoint orbit of $\lambda$.
 More precisely, there is a diffeomorphism
$\O\to  G/P_\lambda$ with the property that the pull-back of the 
Fubini-Study symplectic form on $G/P_\lambda$ (see above) is just Kirillov-Kostant-Souriau
symplectic form on $\O$. We deduce that we have 
$\Psi|_{G/P_\lambda}=\Phi$, where $\Phi$ is given by
$$\Phi: \O \hookrightarrow \k^*\to \t^*,$$
which is the composition of the inclusion map with the restriction map.
The details of this discussion can be found for instance in
  \cite[Chapter V]{Huc}. 
  
The identification (\ref{ide}) leads us to the following approach:
the coadjoint orbit $\O$   admits an action of $G$ and the Schubert cells are the orbits 
of the induced $B$ action. To be more precise, we consider the natural action of the Weyl group 
 $W=N_K(T)/T$ on $\t^*$ and the orbit $W\lambda=W/W_\lambda$, which is contained in $\O$. Then the Schubert cells are $Bw\lambda$
and the Schubert varieties are
$$X(w)=\overline{Bw\lambda},$$ for all $w\in W$ (the closure is relative to the
standard topology of $\O$).
Thus, what we actually study in this paper are symplectic quotients of $X(w)$ relative to the
$T$ moment map  $\Phi:\O\to \t^*$.

 We also consider a circle $S\subset T$ and the moment map
of the $S$ action on $\O$: this is just the composition of  $\Phi$ with 
the restriction map $\t^* \to {\rm Lie}(S)^*$.
More precisely, let us consider $a\in \t$ such that   $S=\exp (\bR a)$.
In fact, $a$ is
an element of the integral lattice of $T$, that is,
$\ker(\exp:\t\to T)$.  We can also assume that $a$  is not an integer multiple of any other integral element.  Denote by $\nu$ the element of the dual space 
$ (\bR a)^*$ determined by $\nu(a)=1$. Then the $S$ moment map of $\O$ is
$ \Phi_a \nu : \O \to (\bR a)^*=\bR \nu$, where  $\Phi_a$ denotes the evaluation of
$\Phi$ at $a$. 

Let us fix a Schubert variety $X=X(w)$, where $w\in W$.
We study the topology of the  quotients 
$$(\Phi^{-1}(\mu)\cap X)/T \  \ \ {\rm and } \ \ \  (\Phi_a^{-1}(r_0)\cap X)/S,$$
 where $\mu\in \Phi(X)$ and $r_0\in \Phi_a(X)$. In the case where $X$ is non-singular,
 these spaces are symplectic quotients of $X$. If $X$ does have singularities,
 we can choose $\mu$ and $r_0$ such that   $\Phi^{-1}(\mu)$ and
 $\Phi_a^{-1}(r_0)$ are contained in $X\setminus {\rm Sing}(X)$. The two quotients
 are symplectic quotients of the latter space, which is a (non-compact) K\"ahler manifold. The assumption above will {\it not} be in force everywhere in this paper.

\noindent{\bf Remark.}  Such quotients could be  relevant for the study of the Demazure ($B$-)module
${ V}_{w}(\lambda)$ (cf. \cite{Br-Ku}).  Guillemin and
Sternberg used ``quantization commutes with reduction" to prove the following
result: if a weight vector $\mu$ occurs among the weights of the $K$ representation
$V_\lambda$ (irreducible of highest weight $\lambda$, see above), then $\mu$ is contained in $\Phi(\O)$; moreover,
the multiplicity of $\mu$ is equal to the dimension of the space of sections of
the line bundle induced by $L_\lambda$ on the symplectic quotient $\Phi^{-1}(\mu)/T$
(cf. \cite{Gu-St}, \cite{Sj}). Now the Demazure module 
${ V}_w(\lambda)$ is the dual of the  space of sections $H^0(X,L_\lambda|_X)$
with the canonical $B$ action
(see \cite[Corollary 3.3.11]{Br-Ku}). We can expect again that if the weight vector
$\mu$ is a weight of $V_w(\lambda)$, then $\mu$ is in $\Phi(X)$ (see also Remark 1 below) and  the multiplicity of $\mu$ is equal to the dimension of the space of sections of
the line bundle induced by $L_\lambda$ on  $(\Phi^{-1}(\mu)\cap X)/T$.
This is certainly not an obvious result. First, because $X$ is not smooth:
however, we could use Teleman's  \cite{Te} ``quantization commutes with reduction" theorem,
which holds for linear group actions
on projective varieties which have rational singularities (Schubert varieties
do have this property). Second, the Borel subgroup $B$ is not reductive: Ion \cite{Ion} was able to overcome this  and obtain geometric formulas for the multiplicities of the weights of $V_w(\lambda)$ by extending methods
which had been used previously by
Mirkovi\'{c} and Vilonen in the context of  representations of reductive Lie groups.
We will not explore such phenomena in this paper.

Our first result states (or rather implies) that the  quotients defined above are connected. Note that the result holds without any assumption on $\mu$ or $r_0$.

\begin{theorem}\label{zeromain} All preimages $\Phi^{-1}(\mu)\cap X$, $\mu\in \Phi(X)$, and  
$\Phi_a^{-1}(r_0)\cap X$,  $r_0\in \Phi_a(X)$, are connected.
\end{theorem}

\noindent {\bf Remarks.} 
1. This theorem is related to the convexity theorem for Hamiltonian torus actions
on symplectic manifolds of Atiyah and Guillemin-Sternberg.
Namely, Atiyah's proof  of the latter result uses the fact that all preimages
of the moment map are connected (cf. \cite{At}, see also   \cite[Section 5.5]{Mc-Sa}).
In the same paper \cite{At}, he shows that  if $X$ and $\Phi$ are as above,  then $\Phi(X)$ is the convex hull  of  the set $\Phi(X^T)$
in $\t^*$.  Since $X$ is in general not smooth, the argument involving the connectivity of the preimages cannot be used: 
instead, Atiyah uses a  convexity result for closures of 
$T^{\bC}$ orbits on K\"ahler manifolds. It would be interesting to find a  proof of the convexity of $\Phi(X)$  which uses Theorem \ref{zeromain}.

2. Here is a simpler proof of the connectivity of $\Phi^{-1}(\mu)$, under the assumption that $\mu$ is a weight vector in $\t^*$. Let us  consider
 the character $\chi_{\mu}: T^{\bC}\to \bC^*$ induced by $\mu$ and the
twisted   action  of $T^{\bC}$ on the line bundle $L_{\lambda}$ (see above), which is defined as follows: $$t.[(g,z)]:=[(gt^{-1},\chi_{\mu}(t)^{-1}z)],$$
for all $t\in T^{\bC}, g\in G$, and $z\in \bC$. Let $X^{ss}(L_{\lambda})/T^{\bC}$ be the 
corresponding  Geometric Invariant Theory (shortly GIT) quotient (here $X^{ss}( L_{\lambda})$ is the space of semistable points of the pair 
$(X,L_{\lambda})$). It is  homeomorphic to the symplectic quotient   $(\Phi^{-1}(\mu)\cap X)/T$, by a theorem of Kirwan and Ness (see e.g. \cite[Section 1]{Hu} and
\cite[Section 2]{Hu1})). 
The GIT quotient is connected, because $X^{ss}(L_{\lambda})$ is connected (in turn, this follows from the fact that $X$ is an irreducible projective variety and $X^{ss}(L_{\lambda})$ is Zariski open in $X$). Thus  
$(\Phi^{-1}(\mu)\cap X)/T$ is connected. Consequently,   $\Phi^{-1}(\mu)\cap X$ is connected as well.
Note that   this proof  works only  in the case where    $\mu$  is
a weight vector. We wanted to mention it because the proof of Theorem \ref{zeromain}  (see Section 2 below) does use these ingredients, together with some others, 
most importantly a theorem of Heinzner and Migliorini \cite{He-Mi}.  
Everything we said here remains true if we replace $T$ by $S$ and $\mu$ by an integer number $n$.

The cohomology of GIT quotients of {\it smooth} projective varieties has been extensively investigated during the past two decades, starting with the seminal work of
Kirwan \cite{Ki1}. By contrary,  little seems to be known in this respect about quotients
of {\it singular} varieties equipped with algebraic group actions. Except the 
results of \cite{Ki2} (where the {\it intersection} cohomology is discussed), we are not aware of any other approaches concerning this topic.  
The next two theorems give a description of the {\it ordinary} (i.e. singular) 
cohomology\footnote{All cohomology rings will be with coefficients in $\bR$.} ring
of the symplectic quotients $$X\quott_{\lambda}S(r_0):=(\Phi_a^{-1}(r_0)\cap X)/S.$$   We can only do that  under certain restrictions on $a$, $r_0$, $\lambda$, and  $w$. 

 \noindent {\bf Assumption 1 (concerning $a$).}  {\it The vector $-a$ is in the (interior of the) fundamental Weyl chamber of $\t$. The numbers
 $\Phi_a(v\lambda)$, $v\lambda \in W\lambda$, are any two distinct.
}
 
 One consequence of this is that the fixed point set $\O^S$ is given by
 $$\O^S =\O^T=W\lambda. $$ 
Thus, for our Schubert variety $X=X(w)$ we have
 $$X^S=X^T=W\lambda \cap X.$$
 Assumption 1 also implies that the unstable manifolds of $\Phi_a$ relative to the
 K\"ahler metric on $\O$ are just the Bruhat cells (cf. e.g. \cite{Du-Ko-Va}). 
 This will  allow us to use
 Morse theory for the restriction of   $\Phi_a$ to $X\setminus {\rm Sing}(X)$,
 which is one of the  main tools we will be employing in our proofs:
 for instance, we will show that the critical points of this function  are  in
 $W\lambda$ (see Lemma \ref{cri}).

An important instrument will be the Kirwan map 
$$\kappa:H^*_S(X)\to H^*_S(\Phi_a^{-1}(r_0)\cap X).$$
The domain of this map can be described as follows.
We start with a Goresky-Kottwitz-MacPherson 
\cite{Go-Ko-Ma} type presentation of  the ring $H^*_T(\O)$
(cf. \cite{Go-Ko-Ma}):
the  map $$H^*_T(\O)\to H^*_T(\O^T)=\bigoplus_{v\lambda \in W\lambda}
H^*_T({\rm pt})$$
 induced by the inclusion $\O^T\subset \O$ 
is injective and we know exactly its image 
 (see for instance    \cite[Section 2.3]{Gu-Ho-Za} or the discussion preceding
 Lemma \ref{katwo} below). The action of $T$ on $X$ is equivariantly formal, by \cite[Section 1.2]{Go-Ko-Ma}:
 indeed, any of the
 Bruhat cells in the CW decomposition of $X$ is $T$ invariant. 
Since $X$ is a closed $T^{\bC}$ invariant subvariety of $\O$,
the GKM description of $H^*_T(\O)$ yields  readily  a similar description of the ring $H^*_T(X)$, as the image of the (injective) map
 \begin{equation}\label{restr}\jmath^*: H^*_T(X)\to H^*_T(X^S)=\bigoplus_{v\lambda \in X^S}
H^*_T({\rm pt}).\end{equation}
induced by the inclusion $\jmath : X^S\to X$.
The image of $\jmath^*$ is described in terms of the moment graph of $X$. 
The vertices of this graph are the elements of the set $X^S=W\lambda \cap X$.
If $\gamma \in \Phi^+$ and $v\in W$ such that 
$v\lambda \neq s_\gamma v\lambda$ and both $v\lambda$ and $s_\gamma v\lambda$ are in $X$, then we join the vertices $v\lambda$ and
$s_\gamma v\lambda$  by an edge, which is labeled with  $\gamma$.
Denote by $\Gamma$ the resulting graph.
The image of $\jmath^*$ consists of all ordered sets
 $(p_{v\lambda})_{v\lambda \in X^S}$ 
 where $p_{v\lambda}\in S(\t^*)$, which are {\it admissible } relative to
 $\Gamma$: by this we mean that  if $v\lambda$ and
 $u\lambda$ are joined by an edge with label $\gamma$, then the
 difference $p_{v\lambda}-p_{u\lambda}$ is divisible by $\gamma$
 (cf. e.g. \cite{Gu-Ho-Za}, see also the discussion preceding Lemma \ref{katwo}
 below).
 Now because the action of $T$ on $X$ is equivariantly formal, 
 the map $H^*_T(X)\to H^*_S(X)$ is surjective.  
 The image of the injective map
$$H^*_S(X)\to H^*_S(X^S)=\bigoplus_{v\lambda \in X^S}
H^*_S({\rm pt})$$
is obtained from the image of the  map given by (\ref{restr}) by projecting it via the canonical map
$$H^*_T({\rm pt})={\rm Symm}(\t^*)\to {\rm Symm} ((\bR a)^*) =H^*_S({\rm pt}).$$

If $X$ were smooth and $r_0$ a regular value of $\Phi_a|_X$, we could determine 
the ring $H^*(X\quott_{\lambda}S(r_0))$ as follows:
use  that the action of $S$ on $\Phi_a^{-1}(r_0)$ is locally free, which implies that we have the ring isomorphism
$H^*_S(\Phi_a^{-1}(r_0)\cap X)\simeq H^*(X\quott_{\lambda}S(r_0))$;
then use that  $\kappa$ is surjective  (cf. Kirwan \cite{Ki1}),  and the Tolman-Weitsman \cite{To-We2} description of $\ker \kappa$. In the second part of this  paper we
give examples of non-smooth Schubert varieties for which this program still works.
Namely, they must satisfy the following assumptions.     

\noindent {\bf Assumption 2 (concerning $X$).} {\it The singular set of $X$ consists of one single point, that is, we have
$${\rm Sing}(X)=\{\lambda\}.$$
}

For example, let us consider 
the Grassmannian $G_2(\bC^n)$ of 2-planes in $\bC^n$. Take $p$ an integer such 
that $4 \le p \le n$. The Schubert variety 
$$X=\{V\in G_2(\bC^n) \ : \  \dim (V\cap \bC^2)\ge 1,
\dim (V\cap \bC^{p})\ge 2\}$$
has one singular point  (to prove this, we  use   \cite[Theorem 9.3.1]{Bi-La}).
The case $n=p=4$ will be discussed in detail in the last section of the paper.

We will also need an assumption concerning $r_0$. 
This is expressed in terms of the moment graph of $X$ (see above).
We  note that, by Assumption 1, the points $\lambda$ and $w\lambda$
are the global minimum point, respectively  the global maximum point  of $\Phi_a|_X$. 
Let us remove from the graph $\Gamma$ all vertices $v\lambda$ with
$\Phi_a(v\lambda)<r_0$, as well as any edge with at least one endpoint at such a
vertex.  Denote by $\Gamma_{r_0}$ the resulting graph.

\noindent{\bf Assumption 3 (concerning $r_0$).} 
{\it (i)The number $r_0$ is in $\Phi_a(X)\setminus \Phi_a(X^S)$.

(ii) If the ordered set $(p'_{u\lambda})_{u\lambda \in X^S, \Phi_a(u\lambda)>r_0}$ is admissible
relative to $\Gamma_{r_0}$, then there exists an ordered set
$(p_{v\lambda})_{v\lambda \in X^S}$ which is admissible relative
to $\Gamma$ and $p_{u\lambda}=p'_{u\lambda}$ whenever
$\Phi_a(u\lambda)>r_0$.
}

The second point of this assumption seems to be  hard to verify. We can always find $r_0$  for which this condition is satisfied. This happens for instance when
$r_0$ is ``high enough", such that $\Gamma_{r_0}$
consists of only one point, that is,
$\Gamma_{r_0}=\{w\lambda\}$: an extension $(p_{v\lambda})_{v\lambda \in X^S}$
of the polynomial $p'_{w\lambda}$ is given by
$p_{v\lambda}=p'_{w\lambda}$ for all $v\lambda \in W\lambda$.
However, we will see in the last section an example where there exist
numbers $r_0$
 such that $\Gamma_{r_0}$ has more than one vertex and Assumption 3 (ii) is satisfied.

These assumptions will allow us to  prove the Kirwan surjectivity theorem for $\kappa$.

\begin{theorem}\label{firstmain} {\rm (Kirwan surjectivity)} If Assumptions 1, 2, and 3  are satisfied, then
the map $\kappa$ is surjective.
\end{theorem}

To determine the cohomology of our quotient, we first
notice that, by  Assumption 3, point (i), $r_0$ is a regular value of the map
$\Phi_a$ restricted to $X\setminus \{\lambda\}$ (see also Lemma \ref{cri} below).
This space is a K\"ahler $S$ invariant  submanifold of $\O$ and $\Phi_a|_{X\setminus \{\lambda\}}$ is a moment map. Thus the action of $S$ on the level $\Phi_a^{-1}(r_0)\cap X$ is locally free and we have
$$H^*_S(\Phi_a^{-1}(r_0)\cap X)\simeq H^*(X\quott_\lambda S(r_0)).$$ 
We also deduce that our quotient $X\quott_\lambda S(r_0)$ has at most orbifold singularities.
 A complete description of  the ring $H^*(X\quott_\lambda S(r_0))$ will be obtained
 after finding the kernel of $\kappa$.  This is done by
the following theorem.

\begin{theorem}\label{thirdmain} {\rm (The Tolman-Weitsman kernel)}   If Assumptions 1, 2,   and 3  are satisfied, then the 
kernel of $\kappa$ is equal to $K_-+K_+$. Here
$K_-$ consists of all $\alpha \in H^*_S(X)$ such that
$$\alpha|_{v\lambda}=0,  \ {\it for \ all} \  v\lambda \in X^S
\ {\it  with} \ 
\Phi_{a}(v\lambda)<r_0,$$
and $K_+$ is defined similarly (the last condition is $\Phi_{a}(v\lambda)>r_0$).
\end{theorem}

Here, by $\alpha|_{v\lambda}$ we have denoted the image of $\alpha$
under the map $H^*_S(X)\to H^*_S(\{v\lambda\})$ induced by the
inclusion $\{v\lambda\}\to X$.

\noindent {\bf Remark.}   
The space $X\quott_\lambda S(r_0)$ is a symplectic quotient of the K\"ahler manifold $X\setminus \{\lambda\}$ with respect to the
$S$ action. However, Theorems \ref{firstmain} and \ref{thirdmain} are not
direct consequences of the classical results known in this context.
For example,   
 the surjectivity of $\kappa$ is not  a consequence of Kirwan's surjectivity theorem
 (cf. \cite{Ki1}, \cite{To-We2}): indeed,  the restriction of the moment map 
$\Phi_a$ to $X\setminus \{\lambda\}$ is not a proper map (because it is bounded, as 
$\Phi_a(\O)$ is a bounded subset of $\bR$). Thus we cannot use Morse theory
for $(\Phi_a-r_0)^2$ on   $X\setminus \{\lambda\}$. 
We use instead  the restriction of the function $\Phi_a$  and its properties, like the fact
that its critical set is $X^S\setminus\{\lambda\}$ and its unstable manifolds are  Bruhat cells
(see Section 3 below). To understand exactly why is  each of the three assumptions necessary for our development, one can see the remark at the end of Section \ref{section3}.

 \vspace{0.2cm}

\noindent {\bf Acknowledgements.} I wanted to thank Lisa Jeffrey for discussions concerning the topics of the paper. I am also grateful to the referee for carefully reading the manuscript and 
suggesting many improvements. 

\section{Connectivity of the levels of $\Phi$}

In this section we will be concerned with Theorem \ref{zeromain}. 
We mention that a similar connectivity result (for Schubert varieties in  loop groups) has been proved in  \cite[Section 3]{Ha-Ho-Je-Ma},  by using essentially the same arguments as here. 

We start with some general considerations. Let $Y$ be a K\"ahler manifold acted on holomorphically by a compact torus $T$. Assume that $T$ preserves the K\"ahler structure  and the action of $T$ is Hamiltonian. Let 
$\Psi: Y\to {\rm Lie}(T)^*$ be any moment map. A point $y\in Y$ is called $\Psi$-semistable if the intersection
$\overline{T^{\bC}y}\cap \Psi^{-1}(0)$ is non-empty, where $T^{\bC}$ is the complexification of $T$ and $\overline{T^{\bC}y}$ is the closure of the orbit of $y$.  We denote by $Y^{ss}(\Psi)$ the space of all $\Psi$-semistable points of $Y$. 
We also choose an inner product on $\t$, denote by
$\| \cdot \|$ the corresponding norm, and consider the function $\|\Psi\|^2 :Y\to \bR$. 
We denote by 
$Y^{\min}(\|\Psi\|^2)$  the minimum stratum of $\|\Psi\|^2$, that is, the space of all points $y\in Y$ with the property that the $\omega$-limit of the integral curve through 
$y$ of the negative gradient vector field $-{\rm grad} (\|\Psi\|^2)$  is contained in
$\Psi^{-1}(0)$.   The following result is a direct consequence of
\cite[Theorem 6.18]{Ki1}.

\begin{theorem}\label{khm} {\rm (Kirwan)} We have  $$Y^{ss}(\Psi)=Y^{\min}(\|\Psi\|^2).$$ 
\end{theorem}

There is another version of the notion of semistability, which is defined as follows.
Let us assume that $Y$ is a smooth projective variety.
We endow $Y$ with the K\"ahler structure induced by its projective embedding.
Let $L$ be a $T^{\bC}$ equivariant  ample line bundle on $Y$.
 A point $y\in Y$ is called $L$-semistable if there exists an integer number $n\ge 1$  and a $T^{\bC}$ equivariant section $s$ of $L^{\otimes n}$ such that $s(y)\neq 0$.
 We denote by $Y^{ss}(L)$ the set of all $L$-semistable points in $Y$.
 The following theorem has been proved by Heinzner and Migliorini \cite{He-Mi}.
 
 \begin{theorem}\label{hm} {\rm (Heinzner and Migliorini)}
  If $Y$ is a smooth projective variety with a holomorphic action of $T$
  which preserves the K\"ahler form
  and $\Psi:Y\to \t^*$ is a moment map, then there exists a very ample $T^{\bC}$  equivariant line bundle $L$ on $Y$ such that 
$$Y^{ss}(\Psi)=Y^{ss}(L).$$
\end{theorem}

The two theorems above will allow us to prove Theorem \ref{zeromain}, as follows. 

\noindent {\it Proof of Theorem \ref{zeromain}.}  We take $\mu\in \Phi(X)$ and show that $\Phi^{-1}(\mu)\cap X$ is connected. To this end, we consider the function
$g:\O \to \bR$, $g(x) = \|\Phi(x) - \mu\|^2$.  By  Theorem \ref{khm}, we have
$$\O^{\min}(g)=\O^{ss}(\Phi-\mu).$$  
Since the action of $T$ on $\P(V_\lambda)$ is linear, it is holomorphic and it leaves the
Fubini-Study symplectic form invariant. By Theorem \ref{hm}, there exists a $T^{\bC}$-equivariant very ample line bundle $L$ on $\O$ such that 
$$\O^{ss}(\Phi-\mu) = \O^{ss}(L).$$ The semistable set of 
$X\subset \O$ with respect to the line bundle $L|_{X}$  is  
$$X^{ss}(L|_X) = \O^{ss}(L)\cap X.$$
This is a Zariski open subspace of $X$. 
The Schubert variety $X=\overline{Bw\lambda}$ is  irreducible: indeed, since
$B$ is connected,  
the orbit $Bw\lambda$ is an irreducible locally closed projective variety.
Consequently, $X^{ss}(L|_X)$ is a connected topological subspace of 
$X$ relative to the differential topology (by
\cite[Corollary 4.16]{Mu}). 
We will need the following claim.

\noindent {\it Claim.}
The space
$X^{ss}(L|_X)=\O^{\min}(g)\cap X$ contains $g^{-1}(0)\cap X$ as a deformation retract. 

To prove this, 
let us first consider the flow $\sigma_t$, $t\in \bR$, on $\O$ induced by the vector
field $-{\rm grad}g$ (the gradient is relative to the K\"ahler metric on $\O$). By a theorem of Duistermaat
(see for instance \cite[Theorem 1.1]{Le}), 
the map $$[0,\infty]\times 
\O^{\min}(g)\to \O^{\min}(g), \ (t,x)\mapsto \sigma_t(x)$$ is a deformation retract
of $\O^{\min}(g)$ to $g^{-1}(0)$. 
The claim follows from the fact that for any $t\in [0,\infty)$, the automorphism
$\sigma_t$ of $\O$    leaves $Bw\lambda$, hence also $X=X(w)=\overline{Bw\lambda}$, invariant. 
Indeed,  for any $x\in Bw\lambda$, we have 
\begin{equation}\label{rmgr}({\rm grad} g)_x=2J_x((\Phi(x)-\mu).x),\end{equation} where $J_x$ denotes the complex structure of $\O$ at $x$.
In the equation above, we use the inner product on $\t$ to identify  $\Phi(x)-\mu$ 
with an element of $\t$; this induces the infinitesimal tangent vector 
$(\Phi(x)-\mu).x$ at $x$ (cf. \cite[Lemma 6.6.]{Ki1}).
Equation (\ref{rmgr}) implies that the vector
$({\rm grad} g)_x$  is tangent to $Bw\lambda$, as this space  is a complex $T^{\bC}$  invariant submanifold of $\O$. The claim is proved.

The claim implies  that the space
$g^{-1}(0)\cap X =\Phi^{-1}(\mu)\cap X$ is connected. 

The fact that $\Phi_a^{-1}(r_0)\cap X$ is connected can be proved in a similar way.
\hfill $\square$

\section{Kirwan surjectivity and the Tolman-Weitsman kernel via Morse theory}\label{section3}
Assumptions 1, 2, and 3 are in force  throughout this section. 
We will prove  Theorems \ref{firstmain} and \ref{thirdmain}.
To this end, let us consider the function $$f:=\Phi_a|_X:X \to \bR.$$
The main instrument of our proofs will be Morse theory for the function $f$ 
restricted to $X\setminus  \{\lambda\}$. This space is a smooth, non-compact
submanifold of the orbit $\O$. Its tangent space at any of its points 
is a complex vector subspace of the tangent space to $\O=G/P_\lambda$ at that point
(since $X$ is a complex subvariety of $G/P_\lambda$). 
Thus the restriction of the canonical K\"ahler structure of $\O$ (cf. e.g. 
\cite[Section 4]{At}) to
 $X\setminus \{\lambda\}$ makes the latter space into a K\"ahler manifold.
We denote by $\langle \ , \ \rangle$ the corresponding Riemannian metric. We start with the following lemma, which is a consequence of Assumption 1.

\begin{lemma}\label{cri} The critical set of $f$ restricted to $X\setminus \{\lambda\}$
is $W\lambda\cap X\setminus \{\lambda\}$. All critical points are non-degenerate.
\end{lemma}

\begin{proof} Let $\psi_t:\O \to \O$, $t\in \bR$, denote the flow on $\O$ determined
by the gradient vector field ${\rm grad} \Phi_a$ with respect to the K\"ahler metric. 
The fixed points of this flow are the critical points of $\Phi_a$, that is, the elements of
$W\lambda$. If $v\lambda$ is such a point, we consider the unstable manifold
$\{x\in \O \ | \ \lim_{t\to \infty}\psi_t(x)=v\lambda\}$. This is the same as the
Bruhat cell $Bv\lambda$ (see \cite[Section 4]{At}). 
Consequently, for any $t\in \bR$, the automorphism $\psi_t$ 
of $\O$ leaves each Bruhat cell invariant. Thus it leaves $X\setminus \{\lambda\}$
invariant (since this space is a union of Bruhat cells). We deduce that
the vector field  ${\rm grad} \Phi_a$ is tangent to  $X\setminus  \{\lambda\}$ at any
of its points. Its value at any such point $x$  must be the same as $({\rm grad}f)_x$.
In conclusion, the critical points mentioned in the lemma are those points $x\in X\setminus \{\lambda\}$
with the property that $({\rm grad} \Phi_a)_x=0$. This condition is equivalent to
$x\in W\lambda$. 

The last assertion in the lemma follows from the fact that $f$ is a moment map
of the $S$ action on $X\setminus  \{\lambda\}$. Thus, it is a Morse function
(cf. e.g. \cite[p. 39]{Ki1}).
\end{proof}

For any number $r$ we  denote 
$$X_r^-:=f^{-1}((-\infty,r)),\quad X_r^+:=f^{-1}((r,\infty)).$$
Since the function $\Phi_a$ and the subspace $X$ of $\O$ are $S$ invariant,
$X^-$ and $X^+$ are $S$ invariant subspaces of $X$. If $\alpha \in H^*_S(X)$ and
$A$ is an $S$ invariant subspace of $X$, we denote by $\alpha|_A$ the image of $\alpha$ under the map $H^*_S(X)\to H^*_S(A)$ induced by the inclusion $A\to X$.
We are now ready to state our next lemma. 

\begin{lemma}\label{lemmafirst} Take
$\epsilon >0$ such that 
$f(\lambda)
<r_0-\epsilon $ and the intersection
$f^{-1}([r_0-\epsilon,r_0+\epsilon])\cap X^S$ is empty.  
 Then we have 
 \begin{equation}\label{kerka}\ker\kappa =K'_-+K'_+\end{equation} where we have denoted
$$K'_-=\{\alpha \in H^*_S(X) \ : \ \alpha|_{X_{r_0+\epsilon/3}^-}=0\},$$
and $$K'_+=\{\alpha  \in H^*_S(X) \ : \ 
\alpha |_{X_{r_0-\epsilon/3}^+}=0\}.$$
\end{lemma}
\begin{proof} To simplify notations, put $X^-:=X_{r_0+\frac{\epsilon}{3}}^-$, 
$X^+:=X_{r_0-\frac{\epsilon}{3}}^+$. 
Both $K_-'$ and $K_+'$, hence also their sum, are evidently contained
in $\ker\kappa$. 
We  show that $\ker\kappa\subset K_-'+K_+'$.
To this end, we need the following claim.

{\it Claim.} The space $X^-\cap X^+=f^{-1}((r_0-\frac{\epsilon}{3}, r_0+
\frac{\epsilon}{3}))$ contains $f^{-1}(r_0)$ as an 
$S$-equivariant  deformation retract.

The idea of the proof is to deform $f^{-1}((r_0-\frac{\epsilon}{3},r_0+
\frac{\epsilon}{3}))$ onto $f^{-1}(r_0)$
in $X\setminus  \{\lambda\}$ along the gradient lines of the  function 
$f|_{X\setminus  \{\lambda\}}$. This is possible since  the preimage
$f^{-1}([r_0-\epsilon,r_0+\epsilon])\cap X\setminus \{\lambda\}=
f^{-1}([r_0-\epsilon,r_0+\epsilon])
$ is compact and does not contain any critical points of $f$ (by Lemma \ref{cri}).
The arguments we will employ in what follows are standard 
(see for instance \cite[Proof of Theorem 3.1]{Mi}).
Here are the details of the construction.
We start with a smooth function $F:\bR \to \bR$ such that:
\begin{itemize}
\item $F(s)=1$ for $s\in (r_0-\frac{2\epsilon}{3},r_0+\frac{2\epsilon}{3})$
\item $F(s)=0$ for $s$ outside the interval 
$[r_0-\epsilon,r_0+\epsilon]$
\end{itemize}
We then consider the function $\rho:X\setminus \{\lambda\}\to \bR$
given by
$\rho(x)=\frac{F(f(x))}{\|({\rm grad}f)_x\|^2}$,
for all $x\in X\setminus \{\lambda\}$ (here $\| \cdot \|$ is the norm induced by the
K\"ahler metric).
The vector field $\rho {\rm grad} f$ on $X\setminus  \{\lambda\}$ vanishes outside
the compact set $f^{-1}([r_0-\epsilon,r_0+\epsilon])$.
By \cite[Lemma 2.4]{Mi}, it generates a flow $\phi_t$, $t\in \bR$, on
 $X\setminus  \{\lambda\}$.  
For any $x\in X\setminus \{\lambda\}$ and any
$t\in \bR$ we have
$$\frac{d}{dt}f(\phi_t(x))=
\langle ({\rm grad}f)_{\phi_t(x)}, \rho(\phi_t(x))({\rm grad}f)_{\phi_t(x)}\rangle 
=F(f(\phi_t(x))).$$
Assume that $x\in f^{-1}( (r_0-\frac{\epsilon}{3},r_0+\frac{\epsilon}{3}))$.
Then we have $$f(\phi_t(x))=t+f(x)$$
for all $t\in (-\frac{\epsilon}{3},\frac{\epsilon}{3})$ (the reason is that both sides
of the equation represent solutions of the same  initial value problem).
Then $$R_\tau(x):=\phi_{\tau(r_0-f(x))}(x),$$
$\tau\in [0,1]$, $x\in f^{-1}((r_0-\frac{\epsilon}{3},r_0+\frac{\epsilon}{3}))$, defines a deformation retract of
$f^{-1}((r_0-\frac{\epsilon}{3},r_0+\frac{\epsilon}{3}))$ onto $f^{-1}(r_0)$.
It only remains to show that for any $\tau\in [0,1]$, the map $R_\tau$ from
$f^{-1}((r_0-\frac{\epsilon}{3},r_0+\frac{\epsilon}{3}))$ to itself is $S$-equivariant.
This follows from the fact that 
for any $t\in \bR$, the automorphism $\phi_t$ of $X\setminus \{\lambda\}$ is
$S$-equivariant. Indeed, the function $f:X\setminus \{\lambda\}\to \bR$ is
$S$-invariant; since $S$ acts isometrically on $\O$, the vector field 
${\rm grad}f$, hence also $\rho{\rm grad}f$, is $S$ equivariant.
  The claim is now completely proved.

The claim implies that the pair
$(X^-\cap X^+,f^{-1}(r_0))$ is $S$-equivariantly homotopy equivalent to $(f^{-1}(r_0),
f^{-1}(r_0))$. Hence we have $H^*_S(X^-\cap X^+,f^{-1}(r_0))=\{0\}$.
 From the long exact sequence of the triple
$(X, X^-\cap X^+, f^{-1}(r_0))$ we deduce 
 that the canonical map
$$ \psi: H^*_S(X,X^-\cap X^+) \to H^*_S(X, f^{-1}(r_0))$$ is 
an isomorphism.

Let us now focus on the proof of  the inclusion
$\ker\kappa\subset K_-'+K_+'$. 
Take $\alpha\in \ker \kappa$, that is, $\alpha\in H^*_S(X)$ such that $\alpha|_{f^{-1}(r_0)}=0$. From the long exact sequence of the pair $(X,f^{-1}(r_0))$ we deduce that there exists $\beta \in H^*_S(X,f^{-1}(r_0))$ whose image via $H^*_S(X,f^{-1}(r_0))\to H^*_S(X)$ is $\alpha$. We set $$\eta=\psi^{-1}(\beta).$$
We use the 
relative Mayer-Vietoris sequence  of the triple $(X,X^-,X^+)$.
Let $i_-^*:
H^*_S(X,X^-)\to H^*_S(X, X^-\cap X^+)$ and $i_+^*:
H^*_S(X,X^+)\to H^*_S(X, X^-\cap X^+)$ be the maps induced by the
obvious inclusions.
Because $X^-\cup X^+=X$, the exactness of the Mayer-Vietoris sequence implies that the
 map $$H^*_S(X, X^-)\oplus H^*_S(X,X^+) \to H^*_S(X, X^-\cap X^+)$$
defined by $$(\eta_1,\eta_2)\mapsto i_-^*(\eta_1)-i_+^*(\eta_2)$$
is an isomorphism. This map is in the top of the following
commutative diagram. 

 \begin{equation}\label{dia1} \begin{array}{c}
\xymatrix{
 H^*_S(X, X^-)\oplus H^*_S(X,X^+)  \ar[r]\ar[d] & H^*_S(X, X^-\cap X^+)
\ar[d]
\\
 H^*_S(X)\oplus H^*_S(X)\ar[r]&
H^*_S(X)}
\end{array}
\end{equation}

There exists $(\eta_1,\eta_2)\in  H^*_S(X, X^-)\oplus H^*_S(X,X^+)$
such that 
$$ i_-^*(\eta_1)-i_+^*(\eta_2)=\eta.$$ The image of $\eta$ via the right-hand side map in the diagram is $\alpha$. Let $(\alpha_1,\alpha_2)$ be the image of 
$(\eta_1,\eta_2)$ via the left-hand side map in the diagram.
The classes $\alpha_1,\alpha_2$ have the property that
$\alpha_1|_{X^-}=0$ and $\alpha_2|_{X^+}=0$. From the commutativity of the diagram
we have $$\alpha=\alpha_1-\alpha_2.$$ This finishes the proof.
\end{proof}

\noindent{\bf Remark.} In the general context of circle actions on compact
symplectic manifolds, 
Tolman and Weitsman  gave a description of the kernel of the Kirwan map
similar to equation (\ref{kerka}) above (see \cite[Theorem 1]{To-We2}).
Their proof is different from the one above.  However, 
they do mention that their theorem can be proved by using the
Mayer-Vietoris sequence of the triple $(X, X^-,X^+)$ 
 (see \cite[Remark 3.5]{To-We2}).
We have used this idea to prove Lemma \ref{lemmafirst} above.

We will characterize $K_+'$ and $K_-'$ separately. In the next lemma we
describe
$K_+'$.

\begin{lemma}\label{rtwo}
(i) For any $r_2>r_1>f(\lambda)$, the space $f^{-1}([r_1,r_2])\cap (X\setminus  \{\lambda\})$ is compact.

(ii)  Take $r\in f(X)$, $r> f(\lambda)$, $r\notin f(X^S)$. Then the restriction map
$$H^*_S(X_r^+)\to H^*_S(X_r^+\cap X^S)$$ is injective and the canonical map
$$H^*_T(X_r^+)\to H^*_S(X_r^+)$$ is surjective.

(iii) We have $K_+'=K_+$.

(iv) The restriction map $$\kappa_1:H^*_S(X_{r_0-\frac{\epsilon}{3}}^+)\to H^*_S(f^{-1}(r_0))$$ is surjective.
\end{lemma}

\begin{proof}  Point (i) follows from the fact that $f^{-1}([r_1,r_2])$ is contained in $X\setminus \{\lambda\}$, thus
$$f^{-1}([r_1,r_2])\cap (X\setminus\{\lambda\})=f^{-1}([r_1,r_2]).$$
The latter space is compact, because it is closed in $X$.

(ii) We use the Morse theoretical arguments of \cite{To-We1}
and \cite{To-We2} (see also \cite[Section 2]{Ha-Ho})
for the  function $f|_{X\setminus \{\lambda\}}$. 
More precisely, let us first note that this function has a maximum at
the point $w\lambda$ (where $w$ is given by $X=X(w)$): this follows from the fact that
$\lambda$ and $w\lambda$ are the minimum, respectively maximum points of 
$f$ on $X$
(cf. \cite[Section 4]{Du-Ko-Va}). To prove the first assertion
 it is sufficient to
 take $r_1,r_2\in \bR\setminus f(X^S)$ such that 
$r<r_1<r_2$ and note that:
\begin{itemize}
\item if  $f^{-1}([r_1,r_2])\cap X^S=\{w\lambda\}$ then
$H^*_S(X_{r_1}^+)\to H^*_S(\{w\lambda\})$ is injective
(since $\{w\lambda\}$ is an $S$-equivariant deformation retract of $X^+_{r_1}$).
\item if $f^{-1}([r_1,r_2])\cap X^S$ is empty, then the map
$H^*_S(X^+_{r_1})\to H^*_S(X^+_{r_2})$ induced by
the inclusion $X^+_{r_2}\to X^+_{r_1}$ is an isomorphism
(since $X^+_{r_2}$ is an $S$-equivariant deformation retract of $X^+_{r_1}$).
\item if the map $H^*_S(X_{r_2}^+)\to H^*_S(X_{r_2}^+\cap X^S)$ is injective
and $f^{-1}([r_1,r_2])\cap X^S=\{v\lambda\}$ for some $v\in W$, then the map
$H^*_S(X_{r_1}^+)\to H^*_S(X_{r_1}^+\cap X^S)$ is injective as well.
\end{itemize}
To prove the last item, let us consider the following commutative diagram:
\begin{equation}\label{exactseq}
\vcenter{\xymatrix{
\cdots
\ar[r] &
H_{S}^*(X^+_{r_1},X^+_{r_2})
\ar[r]^{\textcircled{\small{2}}}\ar[d]^{\simeq} &
H_{S}^*(X^+_{r_1})
\ar[r] \ar[d]^{\textcircled{\small{1}}} &
H_{S}^*(X^+_{r_2})\ar[r] &  \cdots
\\
&
H_{S}^{*-k } (\{v\lambda\})
\ar[r]^{~~  \cup e}  &
H_{S}^*(\{v\lambda\}) &
\\
}}
\end{equation}
Here $k$ is the dimension of the positive space of the Hessian of
$f$ at the point $v\lambda$, call it $T^+_{v\lambda}X$. 
Also, $\simeq$ denotes the  isomorphism obtained by  composing the excision map $H^*_S(X^+_{r_1},X^+_{r_2})\simeq H^*_S(D^k,S^{k-1})$ (where $D^k$, 
$S^{k-1}$ are the unit disk,
respectively unit sphere in $T^+_{v\lambda}X$),  
with the Thom isomorphism $H^*_S(D^k,S^{k-1})\simeq 
 H^{*-k } _S(\{v\lambda\})$. The map $\textcircled{\small{1}}$ is induced by the inclusion of
 $\{v\lambda\}$ in $X^+_{r_1}$. The cohomology class 
  $e\in H_S^{k}(\{v\lambda\})=H^{k}(BS)$ is the $S$ equivariant Euler class of
  $T^+_{v\lambda}X$. 
   Let us note that the group $S$ acts linearly without fixed points on the latter space
  (since $v\lambda$ is an isolated fixed point of the $S$ action).
By the Atiyah-Bott lemma (cf. \cite{At-Bo}), $e$ is a non-zero element of $H^*(BS)$.
  Thus, the multiplication by 
$e$ is an injective endomorphism of $H^*(BS)$. We deduce that the map
$\textcircled{\small{2}}$ is injective and consequently  the long exact sequence of the pair 
$(X^+_{r_1},X^+_{r_2})$ splits into short exact sequences of the form
\begin{equation}\label{exac}
0\longrightarrow H_{S}^*(X^+_{r_1},X^+_{r_2})
\longrightarrow
H_{S}^*(X^+_{r_1})
\longrightarrow 
H_{S}^*(X^+_{r_2})\longrightarrow 0.\end{equation}
Let us consider now the following commutative diagram.
$$
\begin{array}{c}
\xymatrix{
0  \ar[r] & H_S^*(X^+_{r_1},X^+_{r_2}) \ar[r]\ar[d]_{\textcircled{\small{3}}}
& H_S^*(X^+_{r_1}) \ar[r]\ar[d]_{\imath_{1}^*} & H_S^*(X^+_{r_2}) 
\ar[r]\ar[d]_{\imath_{2}^*}
& 0 \\
0  \ar[r] & H_S^*(\{v\lambda\})\ar[r] &
H_S^*( X^+_{r_1}\cap X^S) \ar[r] & H_S^*( X^+_{r_2}\cap X^S) \ar[r]
& 0 }
\end{array}
$$
  where we have identified $H_S^*( X^+_{r_1}\cap X^S, X^+_{r_2}\cap X^S)=H^*_S(\{v\lambda\})$.   
 By hypothesis, the map $\imath_{2}^*$ is injective. 
The map $\textcircled{\small{3}}$ is the same as the composition of $\textcircled{\small 1}$ and
$\textcircled{\small 2}$ (see diagram (\ref{exactseq})), thus it is injective as well. By a diagram chase we deduce that $\imath_{1}^*$ is injective.   

  The second assertion is  proved by the same method as before, by induction over the sublevels of $-f|_{X\setminus  \{\lambda\}}$. This time we confine ourselves to
  show that if 
  $r_1,r_2\in \bR\setminus f(X^S)$  satisfy
$r_2>r_1>f(\lambda)$, $f^{-1}([r_1,r_2])\cap X^S=\{v\lambda\}$, and 
  the map $H^*_T(X_{r_2}^+)\to H^*_S(X_{r_2}^+)$ is surjective, then the map 
$H^*_T(X_{r_1}^+)\to H^*_S(X_{r_1}^+)$ is surjective too.
To this end we first note that  both  restriction maps
$H^*_T(X_{r_1}^+)\to H^*_T(X_{r_2}^+)$ and 
$H^*_S(X_{r_1}^+)\to H^*_S(X_{r_2}^+)$ are surjective: we use  the exact sequence (\ref{exac}) and its analogue for
$T$ equivariant cohomology. Let us consider the following commutative diagram 
$$
\vcenter{\xymatrix{
0
\ar[r] &
H_{T}^*(X_{r_1}^+,X_{r_2}^+)
\ar[r]^{}\ar[d]^{\textcircled{\small{4}}} &
H_{T}^*(X_{r_1}^+)
\ar[r] \ar[d]^{\textcircled{\small{5}}} &
H_{T}^*(X_{r_2}^+)\ar[r]\ar[d]^{\textcircled{\small{6}}} &  0
\\
0 \ar[r] & H_{S}^*(X_{r_1}^+,X_{r_2}^+)\ar[r] &
H_{S}^*(X_{r_1}^+) \ar[r] & H_{S}^*(X_{r_2}^+)  \ar[r]
& 0}}
$$  
The map ${\textcircled{\small{4}}}$ is  surjective: indeed, as before, we have 
$H_{T}^*(X_{r_1}^+,X_{r_2}^+)\simeq H_T^{*-k}( \{v\lambda\})$,
and similarly if we replace $T$ by $S$; thus the map ${\textcircled{\small{4}}}$ is just the
canonical (restriction) map ${\rm Symm} ({\rm Lie}(T)^*) \to {\rm Symm} ({\rm Lie}(S)^*)$. 
A  diagram chase shows that if ${\textcircled{\small{6}}}$ is surjective, then 
${\textcircled{\small{5}}}$ is surjective as well. 

Point (iii) is a straightforward consequence of the first assertion of (ii), by taking 
$r=r_0-\frac{\epsilon}{3}$.

(iv) We prove the surjectivity of $\kappa_1$ inductively, along the sublevels of the function $f|_{X_{r_0-\frac{\epsilon}{3}}^+}$. At the first induction step, we note that for any number $r$ such that
$$r_0+\frac{\epsilon}{3}<r<{\rm min }\{f(v\lambda) \ : \  v\lambda \in X^S, f(v\lambda) >r_0\}$$
the space $f^{-1}((r_0-\frac{\epsilon}{3},r))$ contains  $f^{-1}(r_0)$ as an $S$-equivariant deformation retract. This can be proved exactly like the claim in the proof of Lemma \ref{lemmafirst}. 
Then we consider $r_1,r_2\in \bR\setminus f(X^S)$ such that
$r_0<r_1<r_2$ and show that the map
$$H^*_S(f^{-1}((r_0-\frac{\epsilon}{3}, r_2)))\to H^*_S(f^{-1}((r_0-\frac{\epsilon}{3}, r_1)))$$
is surjective, in each of the following two situations:
\begin{itemize}
\item the intersection $f^{-1}([r_1,r_2])\cap X^S$ is empty.
\item the intersection $f^{-1}([r_1,r_2])\cap X^S$ consists of exactly one point, say
$v\lambda$, where $v\in W$.
\end{itemize}
We use the argument exposed in the proof of point (ii) above.
In the second situation (when $f^{-1}([r_1,r_2])\cap X^S=\{v\lambda\}$) we use the analogue of the short exact sequence  (\ref{exac}).
\end{proof}

The following lemma is the final step towards the proof of Theorem \ref{firstmain}.
Assumption 3,  (ii) is essential in the proof of this lemma.
We will use  the $1$-skeleton of the $T$ action on $\O$.
By definition  (cf. \cite{To-We1}), this consists of all points in $\O$ whose stabilizer have codimension at most 1.
In the case at hand one can describe the 1-skeleton as follows. For any $v\in W$ and any 
$\gamma\in R^+$ such that
$s_\gamma v\lambda \neq v\lambda$, there exists
a subspace ${\mathcal S}_\gamma(v\lambda)\subset \O$,
which is a  metric sphere in
Euclidean space $\k^*$ relative to an ${\rm Ad}^*(K)$ invariant inner product on the latter space. Moreover, it contains $v\lambda$ and $s_\gamma v\lambda$ as antipodal points.
 The 1-skeleton is the union of all these spheres. The torus $T$ leaves 
 ${\mathcal S}_\gamma(v\lambda)$ invariant; in fact, the sphere  is left pointwise fixed by
 the  kernel of  the character $T\to S^1$  induced by $\gamma$ 
 (cf. \cite[Section 2.2]{Gu-Ho-Za}). From the 
 Goresky-Kottwitz-MacPherson theorem (cf. \cite{Go-Ko-Ma}, \cite{To-We1}),
 we deduce that the map $H_T^*(\O)\to H_T^*(W\lambda)$ is injective and its
 image consists of all ordered sets $(p_{v\lambda})_{v\lambda \in W\lambda}$
 with the property that
\begin{equation}\label{pv} p_{v\lambda}-p_{s_{\gamma}v\lambda} \ {\rm is  \ divisible \ by \ } \gamma,\end{equation}
 for all $v\in W$ and $\gamma\in R^+$. If both
 $v\lambda$ and $s_\gamma v\lambda$ are in $X$, then the whole
 ${\mathcal S}_\gamma(v\lambda)$  is contained in  $X$. Thus, the $T$ equivariant cohomology ring of $X$ is isomorphic to the ring of all ordered sets $(p_{v\lambda})_{v\lambda \in X^S}$  with the property that   the condition (\ref{pv}) holds
 for all $v\in W$ and $\gamma\in R^+$ such that $v\lambda$ and $s_\gamma v\lambda$
 are in $X$.
Finally, we note that the space ${\mathcal S}_\gamma(v\lambda)$ is an orbit of the
complex subgroup $SL_2(\bC)_\gamma$ of $G$ (cf. \cite[Chapter V, Section 6]{Hu}), thus it is a K\"ahler submanifold of $\O$.

\begin{lemma}\label{katwo}
The map $$\kappa_2: H^*_S(X)\to H^*_S(X_{r_0-\frac{\epsilon}{3}}^+)$$ is surjective.
Thus, $\kappa=\kappa_1\circ \kappa_2$ is surjective (the map $\kappa_1$ was
defined in Lemma \ref{rtwo}, (iv)).
\end{lemma}

\begin{proof}
It is sufficient to prove that the map 
$\kappa_2^T: H^*_T(X)\to H^*_T(X_{r_0-\frac{\epsilon}{3}}^+)$ is surjective.
This is because in the commutative diagram
$$ \begin{array}{c}
\xymatrix{
 H^*_T(X) \ar[r]\ar[d]^{\kappa_2^T} & H^*_S(X)
\ar[d]^{\kappa_2}
\\
 H^*_T(X_{r_0-\frac{\epsilon}{3}}^+)\ar[r]&
H^*_S(X_{r_0-\frac{\epsilon}{3}}^+)}
\end{array}
$$
the horizontal maps are surjective (here we have used that $X$ is $T$ equivariantly formal, respectively  Lemma \ref{rtwo}, (ii)).

 Take $\alpha \in H^*_T(X_{r_0-\frac{\epsilon}{3}}^+)$. Take $\gamma\in R^+$
 and $u\in W$ such that   $u\lambda$ and $s_\gamma u\lambda$ are in
 $X$,  $u\lambda \neq s_{\gamma}u\lambda$. We will need the following claim.
 
 \noindent{\it Claim.} If $ u\lambda$ and
 $s_{\gamma}u\lambda$  are  in $X^+_{r_0-\frac{\epsilon}{3}}$,
 then the whole ${\mathcal S}_\gamma(u\lambda)$ is contained in $X_{r_0-\frac{\epsilon}{3}}^+$. 
 
 Indeed, the sphere ${\mathcal S}_\gamma(u\lambda)$  is contained in 
 $X$. Moreover, the restriction of $\Phi_a$ to ${\mathcal S}_\gamma(u\lambda)$ is just the moment map of the $S$ action. Thus its critical points are the fixed points of the
 $S$ action, namely $u\lambda$ and $s_\gamma u\lambda$. One of these two points is a global minimum point. The claim is now proved.  
 
The class $\alpha$ restricted to ${\mathcal S}_\gamma(u\lambda)$ is an element of its $T$ equivariant cohomology. From the discussion preceding the lemma
we deduce that 
$$\alpha|_{ u\lambda}-\alpha|_{s_{\gamma}u\lambda}  \ {\rm is \ divisible \ by} \   \gamma.$$
Here we have denoted by $\alpha|_{ u\lambda}$ the image of $\alpha$ under the map
induced by the inclusion $\{u\lambda\}\to X^+_{r_0-\frac{\epsilon}{3}}$, and similarly for
$\alpha|_{s_{\gamma}u\lambda}$.
Consequently, the ordered set $(\alpha|_{ u\lambda})_{u\lambda \in X^S,
\Phi_a(\u\lambda) >r_0}$ is admissible relative to $\Gamma_{r_0}$.
By Assumption 3, (ii), there exists an ordered set
 $(p_{v\lambda})_{v\lambda \in X^S}$ 
 which is admissible relative to $\Gamma$, such that
  \begin{equation}\label{pw} p_{v\lambda}=\alpha|_{v\lambda}, \ {\rm  whenever} \  v\lambda \in X_{r_0-\epsilon}^+.\end{equation}
 The collection $(p_{v\lambda})$ represents a cohomology class, call it $\beta$, in $H^*_T(X)$.
 By equation (\ref{pw}) and Lemma \ref{rtwo} (ii), we have $\kappa_2(\beta)=\alpha$.
 This finishes the proof.
 \end{proof}

The only piece of information which is still missing is the fact that   $K_-'$ equals
$K_-$. This is the content of the following lemma.

\begin{lemma}\label{under} We have $K_-'=
K_-$.
\end{lemma}

\begin{proof} We only need to show that $K_-\subset K_-'$. 
We actually show that for any $r>f(\lambda)$, $r\notin f(X^S)$,  the map 
$H^*_S(X^-_r) \to H^*_S(X^-_r\cap X^S)$ is injective.
The idea we will use  is the same as in the proof of Lemma \ref{rtwo}, (ii).  
Some adjustments are necessary, though, since $X_r^-$ and all
the other sublevels involved in the argument contain the singular point
$\lambda$. We proceed by proving the following claims. 

\noindent {\it Claim 1.} Take $\delta >0$ such that 
 $\Phi_a(\lambda)<\delta <\Phi_a(v\lambda)$ for any $v\lambda\in W\lambda\setminus\{\lambda\}$. Then the restriction map
$H^*_S(X_{\delta}^-)\to H^*_S(\{\lambda\})$ is injective.

\noindent {\it Claim 2.} If $r_1,r_2\in \bR\setminus f(X^S)$ such that 
$f(\lambda)<r_1<r_2$ and  $f^{-1}([r_1,r_2])\cap X^S$ is empty, then the map
$H^*_S(X^-_{r_2})\to H^*_S(X^-_{r_1})$ induced by
the inclusion $X^-_{r_1}\to X^-_{r_2}$ is an isomorphism.

\noindent {\it Claim 3.} If $r_1,r_2\in \bR\setminus f(X^S)$ such that 
$f(\lambda)<r_1<r_2$, the map $H^*_S(X_{r_1}^-)\to H^*_S(X_{r_1}^-\cap X^S)$ is injective,
and $f^{-1}([r_1,r_2])\cap X^S=\{v\lambda\}$ for some $v\in W$, then the map
$H^*_S(X_{r_2}^-)\to H^*_S(X_{r_2}^-\cap X^S)$ is injective as well.

To prove Claim 1,  we note that $X^-_{\delta}$ is a subset of the
stable manifold at $\lambda$ of the function $\Phi_a: \O \to \bR$ relative to the K\"ahler metric on $\O$.
It is also invariant under 
any of the automorphisms $\psi_t$ of $\O$, where $t\le 0$ (by
$\psi_t$, $t\in \bR$, we  denote the flow of $\O$ induced by 
${\rm grad}\Phi_a$, like in the proof of Lemma \ref{cri}). 
Consequently, it is contractible to $\{\lambda\}$. Moreover, since
$X_\delta^-$ is $S$ invariant, the retract $X_\delta^-\to \{\lambda\}$ is $S$ equivariant. Thus the inclusion map $\{\lambda\} \to X_\delta^-$ induces an isomorphism
$ H^*_S(X_\delta^-) \stackrel{\simeq}{\to} H^*_S(\{\lambda\})$.
 The claim is now proved.

To prove Claim 2, we take into account that $f^{-1}([r_1,r_2])$ is a compact subset of
the smooth manifold $X\setminus \{\lambda\}$. Consequently, there exists a family of
maps $Q_\tau: X_{r_2}^-\setminus \{\lambda\} \times [0,1] \to X_{r_2}^-\setminus \{\lambda\}$, $\tau\in [0,1]$, 
which is a strong deformation retract of $X_{r_2}^-\setminus \{\lambda\}$ on
$X_{r_1}^-\setminus \{\lambda\} $, that is,  $Q_1(X_{r_2}^-\setminus \{\lambda\})\subset
X_{r_1}^-\setminus \{\lambda\}$, 
 $Q_0$ is the identity map of $X_{r_2}^-\setminus \{\lambda\}$, and
 $Q_\tau$ restricted to $X_{r_1}^-\setminus \{\lambda\}$ is the identity map, for all
 $\tau\in [0,1]$.
The concrete expression of $Q_\tau$ can be found for instance in \cite[p. 13]{Mi}: it can be seen from there that for any $\tau\in [0,1]$, the map
$Q_\tau$ is $S$ equivariant. We  extend   $Q_\tau$ to
a continuous map $X_{r_2}^- \to X_{r_2}^-$ by setting $Q_\tau(\lambda)=\lambda$.
This gives  a (continuous) deformation retract of $X_{r_2}^-$ on
$X_{r_1}^-$, which is $S$  equivariant. In this way we have proved Claim 2. 

We now prove Claim 3. Like for the previous claim,  we use Morse theory for the function $f$ on  $X\setminus \{\lambda\}$. By Lemma \ref{cri}, $v\lambda$ is a 
non-degenerate critical point of this function. Denote by $m$ its index.
There exists a closed  $m$-cell 
$e^m\subset X_{r_2}^-\setminus \{\lambda\}$ which can
be attached to $X_{r_1}^-$, and a strong deformation retract $Q'_\tau$, $\tau\in [0,1]$,
of $X_{r_2}^+\setminus\{\lambda\}$ onto $(X_{r_1}^-\cup e^m)\setminus\{\lambda\}$.
From the exact expression of $e^m$ and $Q'_\tau$ (see 
\cite[Proof of Theorem 3.2]{Mi}), we can see that the former is $S$ invariant and the latter is $S$ equivariant. By setting $Q'_\tau(\lambda)=\lambda$, for all $\tau\in [0,1]$, we obtain an $S$ equivariant deformation retract of $X_{r_2}^-$ onto 
$X_{r_1}^-\cup e^m$. We repeat the argument in the proof of Lemma \ref{rtwo} (ii).
\end{proof}

Theorems \ref{firstmain} and \ref{thirdmain} are now completely proved.

\noindent {\bf Remark.} At this point we can understand exactly where each of the assumptions made in the introduction has been used. Namely, Lemmas \ref{cri},
\ref{lemmafirst}, and \ref{rtwo} use only Assumption 1 and  the regularity hypothesis on $r_0$ given in Assumption 3, (i). Lemma 
\ref{under} needs Assumption 2 (see especially Claim 1 in the proof of this lemma). Lemma \ref{katwo} needs Assumption 3, (ii).

\section{An example}

\begin{figure}[h]
\begin{center}
\epsfig{figure=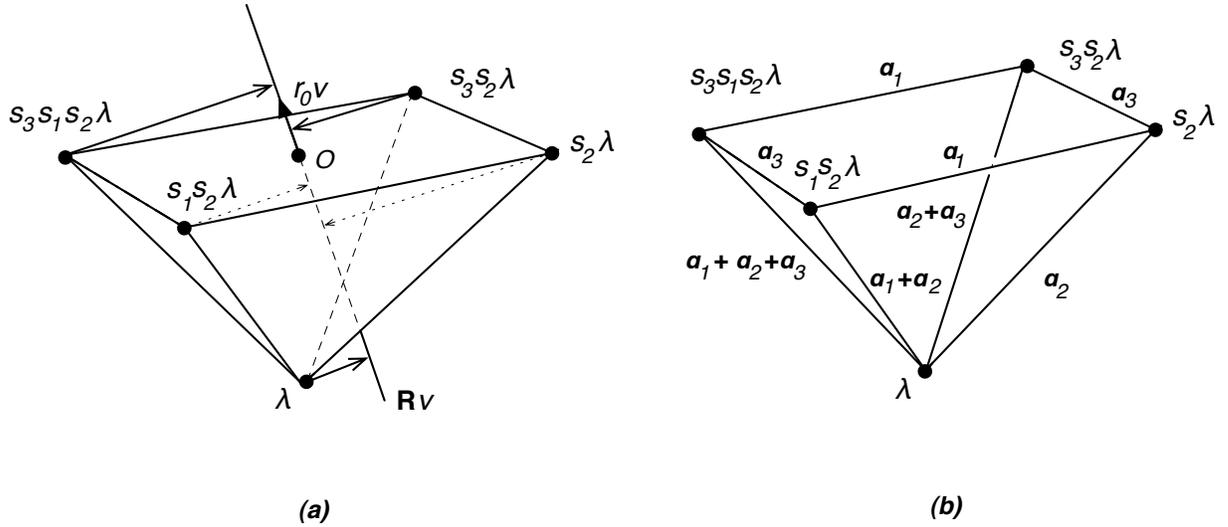,width=6.8in}
\end{center}
\caption{The upside-down pyramid in figure (a) is the image under the moment map of the Schubert variety $X(s_3s_1s_2)$  in the Grassmannian $G_2(\bC^4)$.
Its corners $\lambda$, $s_2\lambda$, $s_1s_2\lambda$, $s_3s_2\lambda$ and $s_3s_1s_2\lambda$ are projected orthogonally onto the line $\bR \nu$, where
$\nu$  is  in the Weyl chamber opposite to the fundamental one. We can also see the vector $r_0\nu$, where
$r_0>0$ satisfies equation (\ref{phid}). Figure (b) describes the moment graph of
$X(s_3s_1s_2)$.}
\end{figure}
 
 The Grassmannian $G_2(\bC^4)$ of 2-dimensional vector subspaces in $\bC^4$ can be identified with a  coadjoint orbit of $SU(4)$. Let us be more specific. 
 Denote by $T$ the space  of all diagonal matrices in $SU(4)$, which is a maximal torus.
Its Lie algebra $\t$ can be described as
 $$\t=\{x=(x_1,x_2,x_3,x_4)\in \bR^4 \ : \ x_1+x_2+x_3+x_4=0\},$$
which is a hyperplane   in $\bR^4$. We equip $\bR^4$, and also its subspace $\t$, with the  Euclidean metric $( \ , \ )$. 
To each $x\in \bR^4$ we assign the element $x^*\in\t^*$ given by 
$$x^*(y):=(x, y),$$
for all $y\in \t$.
A simple root system of $SU(4)$ relative to $T$ is
$$\alpha_1:=\e_2-\e_1,\alpha_2:=\e_3-\e_2, \alpha_3:=\e_4-\e_3,$$
where $e_1,e_2,e_3,e_4$ denotes the standard basis of $\bR^4$.
The fundamental weights are $\e_1, \e_1+\e_2,\e_1+\e_2+\e_3$. The polyhedral cone
generated by them is the fundamental Weyl chamber. 
The coadjoint orbit of any $\lambda$ which is a positive integer multiple of $\e_1+\e_2$ 
can be identified with the Grassmannian $G_2(\bC^4)$
of all 2-planes $V$ in $\bC^4$
(see \cite[Section 5.3]{Gu-Le-St}).
As usual, we denote by $s_i$ the reflection of $\t^*$ corresponding to the root 
$\alpha_i$, where $i\in\{1,2,3\}$. They generate the Weyl group $W$. To each element of $W\lambda$ corresponds a Schubert variety in $G_2(\bC^4)$. In this example we  discuss the Schubert variety $X$ corresponding to $s_3s_1s_2\lambda$. 
This can be  described explicitly as
$$X=\{V\in G_2(\bC^4) \ : \ \dim (V\cap \bC^2)\ge 1\}.$$
It has an isolated singularity at  $\bC^2$, which is the same as 
the point $\lambda$ on the coadjoint orbit (see the paragraph following Assumption 2
in the introduction). So Assumption 2  is  satisfied.  The image of $X$ under the moment map
$\Phi$ is the polytope generated by the elements of $X^T$, which are $\lambda$,  
$s_2\lambda$, $s_1s_2\lambda$,
$s_3s_2\lambda$, and $s_3s_1s_2\lambda$.

We identify $H^*_T({\rm pt})$ with the polynomial ring
$S(\t^*)=\bR[\alpha_1,\alpha_2,\alpha_3]$.
According to the GKM theorem (see the introduction)  the $T$ equivariant cohomology ring $H^*_T(X)$ consists of all elements 
$(f_1,f_2,f_3,f_4,f_5)$ of $H^*_T(\{\lambda\})\oplus H^*_T(\{s_2\lambda\})\oplus H^*_T( \{s_1s_2\lambda\})
\oplus H^*_T(\{s_3s_2\lambda\}) \oplus H^*_T(\{s_3s_1s_2\lambda\})
=S(\t^*)\oplus S(\t^*)\oplus S(\t^*)\oplus S(\t^*)\oplus S(\t^*)$ 
such that
\begin{align}\label{divisi}
{}& f_5-f_4 \ {\rm is \ divisible \ by \ } \alpha_1 \nonumber \\
{}& f_5-f_3 \ {\rm is \ divisible \ by \ } \alpha_3 \nonumber\\
{}&  f_5-f_1 \ {\rm is \ divisible \ by \ } \alpha_1+\alpha_2+\alpha_3 \nonumber\\
{}& f_4-f_2 \ {\rm is \ divisible \ by \ } \alpha_3\\
{}& f_4-f_1 \ {\rm is \ divisible \ by \ } \alpha_2+\alpha_3\nonumber \\
{}& f_3-f_2 \ {\rm is \ divisible \ by \ } \alpha_1\nonumber\\
{}& f_3-f_1 \ {\rm is \ divisible \ by \ } \alpha_1+\alpha_2 \nonumber\\
{}& f_2-f_1 \ {\rm is \ divisible \ by \ } \alpha_2 \nonumber
\end{align}

We  compute the cohomology ring of the symplectic quotient of $X$ with respect to   $S=\exp(\bR a)$, where $a$ is an integral element of $\t$.
The situation is described in Figure 1 (a), where  the identification
$\t^*=\t$ induced by $( \ , \ )$ is in force. We also consider the element $\nu \in \t^*$
such that $\nu(a)=1$ and $\nu$ is identically zero on the orthogonal complement of $a$ in  
$\t$. Then $\nu$ is identified with the vector $a/(a,a)$. For any $x\in\t$, the orthogonal projection of $x^*$ on the line $\bR \nu$ 
is $(x,a)\nu$, which is the same as $\Phi_a(x^*)\nu$. 
The figure shows that    $r_0$ has been chosen such that
\begin{equation}\label{phid}\Phi_a(\lambda) < \Phi_a(s_2\lambda) <\Phi_a(s_1s_2\lambda) <\Phi_a(s_3s_2\lambda)
<r_0<\Phi_a(s_3s_1s_2\lambda).\end{equation}
Assumption 3 is obviously satisfied in this case (see the paragraph following
Assumption 3 in the introduction).
This is the case where we will compute the cohomology ring.
Before doing this, we would like to point out that if $r_0$ satisfies
 $$\Phi_a(\lambda) < \Phi_a(s_2\lambda) <\Phi_a(s_1s_2\lambda) <r_0<\Phi_a(s_3s_2\lambda)
<\Phi_a(s_3s_1s_2\lambda)$$
then Assumption 3 is also satisfied. Indeed, the moment graph is described in Figure 1 (b).
Thus, $\Gamma_{r_0}$ consists of the vertices $s_3s_2\lambda$ and
$s_3s_1s_2\lambda$ which are joined by an edge with label $\alpha_1$.
Let $f'_5$ and $f'_4$ be two polynomials such that
$f'_5-f'_4$ is divisible by $\alpha_1$, that is,
$f'_4=f'_5+\alpha_1 g$, for some $g\in \bR[\alpha_1,\alpha_2,\alpha_3]$.
The polynomials
$$f_5=f'_5, f_4=f'_4=f_5+\alpha_1g, f_3=f_5+\alpha_3g, f_2=f_5+(\alpha_1+\alpha_3)g, f_1=f_5+(\alpha_1+\alpha_2+\alpha_3)g$$ 
satisfy the equations (\ref{divisi}).

Let us now return to the case where $r_0$ satisfies equation (\ref{phid})
and compute the cohomology of $X\quott_{\lambda}S(r_0)$.
To this end we first note that equations (\ref{divisi}) yield the following description of
$H^*_T(X)$: it consists of all elements of  $H^*_T(\{\lambda\})\oplus H^*_T(\{s_2\lambda\})\oplus H^*_T( \{s_1s_2\lambda\})
\oplus H^*_T(\{s_3s_2\lambda\}) \oplus H^*_T(\{s_3s_1s_2\lambda\})$ which are of the form
\begin{align*}{}&(p_1,p_1+\alpha_2p_2,p_1+(\alpha_1+\alpha_2)p_2+\alpha_1(\alpha_1+\alpha_2)p_3,p_1+(\alpha_2+\alpha_3)p_2+\alpha_3(\alpha_2+\alpha_3)p_4,\\{}&
p_1+(\alpha_1+\alpha_2+\alpha_3)p_2+\alpha_1(\alpha_1+\alpha_2+\alpha_3)p_3
+\alpha_3(\alpha_1+\alpha_2+\alpha_3)p_4+\alpha_1\alpha_3(\alpha_1+\alpha_2
+\alpha_3)p_5),\end{align*}
where $p_1,p_2,p_3,p_4,p_5$ are in $\bR[\alpha_1,\alpha_2,\alpha_3]$. 
The restriction map ${\rm Symm}(\t^*)=\bR[\alpha_1,\alpha_2,\alpha_3]\to
  {\rm Symm}((\bR a)^*)=\bR[\nu]$ is given by $\alpha_i\mapsto a_i \nu$, $1\le i \le 3$, where
  $a_i:=\alpha_i(a)$. We deduce that $H^*_S(X)$ consists of the elements of   
  $H^*_S(\{\lambda\})\oplus H^*_S(\{s_2\lambda\})\oplus H^*_S( \{s_1s_2\lambda\})
\oplus H^*_S(\{s_3s_2\lambda\}) \oplus H^*_S(\{s_3s_1s_2\lambda\})$ which are of the form
\begin{align}\label{q1}{}&(q_1,q_1+a_2\nu q_2,q_1+(a_1+a_2)\nu q_2+a_1(a_1+a_2)\nu^2q_3,\nonumber\\
{}&q_1+(a_2+a_3)\nu q_2+a_3(a_2+a_3)\nu^2q_4,\\
{}&q_1+(a_1+a_2+a_3)\nu q_2 +a_1(a_1+a_2+a_3)\nu^2 q_3+a_3(a_1+a_2+a_3)\nu^2
q_4+\nu^3 q_5)\nonumber,\end{align}
where $q_1,q_2,q_3,q_4,q_5\in \bR[\nu]$. 
By Theorem \ref{thirdmain}, the ring $H^*(X\quott_{\lambda}S(r_0))$ is the quotient of $H^*_S(X)$ by the ideal $K_-+K_+$.
Here
$K_-$ is the ideal generated by $(0,0,0,0, \nu^3)$
and $K_+$ consists of all ordered sets of the form (\ref{q1}) where 
$$ q_1+(a_1+a_2+a_3)\nu q_2 +a_1(a_1+a_2+a_3)\nu^2 q_3+a_3(a_1+a_2+a_3)\nu^2q_4+\nu^3 q_5=0.$$

To achieve a better understanding of this ring, we first note that it is graded, with
the graduation  given by $\deg \nu =2$. Thus,
we have $H^{2k+1}(X\quott_{\lambda}S(r_0))=\{0\}$, for all $k\ge 0$.
We also have as follows:
\begin{itemize}
\item The elements of $H^0(X\quott_{\lambda}S(r_0))$ are the
cosets of $(c,c,c,c,c)$, where $c\in \bR$. Thus, $\dim   H^0(X\quott_{\lambda}S(r_0))= 1$, with a basis consisting of the coset of $(1,1,1,1,1)$.
\item The elements of $H^2(X\quott_{\lambda}S(r_0))$ are the cosets of
$$(c_1\nu, c_1\nu+c_2 a_2\nu, c_1\nu+(a_1+a_2)c_2\nu, c_1\nu + (a_2+a_3)c_2\nu,
c_1\nu+(a_1+a_2+a_3)c_2\nu),$$
where $(c_1,c_2)\in \bR^2$. Thus  $H^2(X\quott_{\lambda}S(r_0))$ can be identified
with the quotient of  $\bR^2$ by the kernel of the function
$(c_1,c_2)\mapsto c_1+(a_1+a_2+a_3)c_2$. A basis of this quotient is the
coset of $(c_1,c_2)=(1,0)$.  We deduce that
$\dim   H^2(X\quott_{\lambda}S(r_0))= 1$, with a basis consisting of the coset of $(\nu,\nu,\nu,\nu,\nu)$ in $H^*_S(X)$.
\item Similarly, $\dim   H^4(X\quott_{\lambda}S(r_0))= 1$, with a basis consisting of the coset of $(\nu^2,\nu^2,\nu^2,\nu^2,\nu^2)$.
\item $H^m(X\quott_{\lambda}S(r_0))=\{0\}$ for all $m\ge 5$.
\end{itemize}

From the description above we can see that $X\quott_{\lambda}S(r_0)$ has the same cohomology ring as $\bC P^2$. 
We do not know whether it is actually diffeomorphic to $\bC P^2$ or even whether it is 
smooth.


  
\bibliographystyle{numbering}

\end{document}